\documentclass{amsart}

\usepackage{amssymb}
\newtheorem{defn}{Definition}
\newtheorem{rem}{Remark}

\begin{document}

\title[Embedded Hardy-Littlewood Constants]{Hardy-Littlewood Constants Embedded into Infinite Products over All Positive Integers}

\author{Richard J. Mathar}
\urladdr{http://www.strw.leidenuniv.nl/~mathar}
\email{mathar@strw.leidenuniv.nl}
\address{Leiden Observatory, Leiden University, P.O. Box 9513, 2300 RA Leiden, The Netherlands}


\date{\today}

\begin{abstract}
A group of infinite products over low-order rational polynomials evaluated
at the sequence of prime numbers is loosely called the Hardy-Littlewood constants.
In this manuscript we look at them
as factors embedded in a super-product over primes, semiprimes, 3-almost primes etc.
Numerical tables are derived by transformation into series over $k$-almost prime
zeta functions. Alternative product representations in a basis of
$k$-almost prime products associated with Euler's formula for the
Riemann zeta function
are also pointed out.
\end{abstract}

\maketitle

\section{Appetizer}\label{sec.intro}

As pointed out by J. Vos Post \cite{Vospost112407},
infinite products over rational polynomials
evaluated at integers $n$ can be reordered as the integers are individually
classified as $k$-almost primes \cite{BorweinAMM106}\cite[\S 1.2]{BorweinExpMath},
\begin{equation}
\prod_{n=2}^\infty \frac{n^2-1}{n^2+1}
=
\prod_{k=1}^\infty
\prod_{\substack{n\ge 2\\ \Omega(n)=k}} \frac{n^2-1}{n^2+1}
\label{eq.JVP}
.
\end{equation}
This classification uses
\begin{defn} (Big-Omega, number of prime factors counted with multiplicity)
\begin{equation}
\Omega(n)\equiv \sum_j e_j
\end{equation}
where $n=2^{e_1}3^{e_2}5^{e_3}7^{e_4}11^{e_5}\cdots$ is the prime factorization of $n$.
\end{defn}
One approach to numerical practise is the logarithm of the individual terms,
\cite[1.513.1]{GR}
\begin{gather}
\log \prod_{n=2,\Omega(n)=k}^\infty \frac{n^s-1}{n^s+1}
=
\sum_{n=2,\Omega(n)=k}^\infty \log\frac{n^s-1}{n^s+1}
=
- \sum_{n=2,\Omega(n)=k}^\infty \log\frac{1+1/n^s}{1-1/n^s}
\label{eq.logcycl}
\\
=
-2 \sum_{l=1}^\infty \frac{1}{2l-1}\sum_{n=2,\Omega(n)=k}^\infty \frac{1}{n^{s(2l-1)}}
\nonumber
=
-2 \sum_{l=1}^\infty \frac{1}{2l-1}P_k(s(2l-1)),
\end{gather}
supposed we know how to compute 
the $k$-almost prime zeta functions of
\begin{defn} ($k$-almost-prime Zeta Function, \cite{MatharArxiv0803})
\begin{equation}
P_k(s)\equiv \sum_{\Omega(n)=k} \frac{1}{n^s}.
\label{eq.P1}
\end{equation}
\end{defn}
Examples of this factorization of (\ref{eq.JVP}) are gathered in Table \ref{tab.JVP}.
\begin{table}
\caption{Products emerging from (\ref{eq.JVP}).
More follow
inserting items of Table \ref{tab.zetak}
into (\ref{eq.Zofzeta2}).
}
\begin{tabular}{|l|l|l|}
\hline
$k$ & $s$ & $\prod_{\Omega(n)=k} (n^s-1)/(n^s+1)$ \\
\hline
          1 & 2 &  0.4 = 2/5
\\
          2 & 2 &  0.754499701709514078355718168950541987025077644358722338909979
$\ldots$\\
          3 & 2 &  0.925857274712893127998882138207158415278450218191966021532765
$\ldots$\\
          4 & 2 &  0.980180131878250176004512771774827343256288810830408969770989
$\ldots$\\
\hline
          1 & 3 &  0.704072487320784478296298199978624458092583781119988293242884
$\ldots$\\
          2 & 3 &  0.953501120267507619195724656948780658097471099446350601963684
$\ldots$\\
          3 & 3 &  0.993919830234581988800555498147896650698380465169147914350777
$\ldots$\\
\hline
          1 & 4 &  0.8571428571428571428571428571428571428571428571428571
$\ldots$ = 6/7\\
          2 & 4 &  0.990060339240150340510303763520931737419783647516969651298673
$\ldots$\\
          3 & 4 &  0.999371342684423825202796134016609821126691137699898568119928
$\ldots$\\
\hline
          1 & 5 &  0.930967939958520284033728872221718554673617521060811182585725
$\ldots$\\
          2 & 5 &  0.997730553624228952992631533134776646504723190210043449878780
$\ldots$\\
          3 & 5 &  0.999928848024661518528748225035933938760888742676503394739919
$\ldots$\\
\hline
          1 & 6 &  0.9664335664335664335664335664335664335664335664335664
$\ldots$= 691/715\\
          2 & 6 &  0.999462730036393010562666801735890308203188619330023373330279
$\ldots$\\
\hline
          1 & 7 &  0.983568510511728007936144959029969386664664363392162577192746
$\ldots$\\
          2 & 7 &  0.999870142148595354702394013511164698444556351264029248220202
$\ldots$\\
\hline
          1 & 8 &  0.9919100507335801453448512272041683806389688742629919
$\ldots$ = 7234/7293\\
          2 & 8 &  0.999968223465111892640840105738664881674415326754216354245125
$\ldots$\\
\hline
\end{tabular}
\label{tab.JVP}
\end{table}
The particular case of $k=1$, product over the primes,
results in rational numbers if $s$ is even
(chapter 7 of Ramanujan's first notebook),
\begin{equation}
\prod_p \frac{p^s-1}{p^s+1}
=
\prod_p \frac{1-p^{-s}}{1+p^{-s}}
=
\frac{\zeta(2s)}{\zeta^2(s)}
=
\frac{2B_{2s}}{\binom{2s}{s}B_s^2},
\label{eq.Rama}
\end{equation}
with $B_.$ the Bernoulli numbers, $\zeta(.)$ Riemann's zeta function.
\begin{proof}
This follows from Euler's formula
\begin{equation}
\frac{1}{\zeta(s)}=\prod_p\left(1-\frac{1}{p^s}\right)
,
\label{eq.eul}
\end{equation}
the sign-switched sibling \cite{FrobergBIT8,SebahGourdon}
\begin{equation}
\prod_p \frac{1}{1+p^{-s}}
=
\prod_p \frac{1-p^{-s}}{(1+p^{-s})(1-p^{-s})}
=
\prod_p \frac{1-p^{-s}}{1-p^{-2s}}
=\frac{\zeta(2s)}{\zeta(s)},
\label{eq.zetadoubl}
\end{equation}
and (again Euler's)
\cite[23.2.16]{AS}
\begin{equation}
\zeta(2n)=\frac{(2\pi)^{2n}}{2(2n)!}\left|B_{2n}\right|,\quad n=1,2,3,.\ldots
\end{equation}
\end{proof}
\begin{rem}
The numerical approach to infinite products of rational
polynomials as (\ref{eq.JVP}) is factorization of numerator and denominator
over $\mathbb{C}$ with root sets $\alpha_i$ and $\beta_i$,
use of the limit of $\Gamma$-function ratios \cite[(6.1.47)]{AS}---which
ought evaluate to unity in the sense $\sum (\alpha_i-\beta_i)=0$, easily verified
by subtraction of penultimate coefficients related by the law of Vieta---to end
up with a product over $\Gamma$-ratios \cite[(\S 1.3)]{ErdelyiI}: 
\begin{equation}
\prod_{n=1}^\infty \frac{\prod_i (n-\alpha_i)}{\prod_i (n-\beta_i)}
=
\frac{\prod_i \Gamma(1-\beta_i)}{\prod_i\Gamma(1-\alpha_i)}
\lim_{n\to \infty}
\prod_i
\frac{\Gamma(n-\alpha_i)}{\Gamma(n-\beta_i)}
=
\frac{\prod_i \Gamma(1-\beta_i)}{\prod_i\Gamma(1-\alpha_i)}
.
\label{eq.ntoGamma}
\end{equation}
For lucky roots of the polynomials, further simplification may be possible
by use of the functional equations of the Gamma-function
\cite[(\S 1.2)]{ErdelyiI}\cite{CohenAMS11}.
\end{rem}

\section{Almost-Prime Zeta Functions of the 2nd Kind}
Formulas like (\ref{eq.Rama}) and (\ref{eq.zetadoubl}) indicate that the products defined
in (\ref{eq.JVP}) might not represent the atomic constituents of this arithmetic.
A vague hope to establish some useful arithmetic basis---plus a glimpse
at Euler's formula (\ref{eq.eul})---proposes
\begin{defn}
($k$-almost prime zeta functions of the 2nd kind)
\begin{equation}
\frac{1}{\zeta_k(s)}
\equiv
\prod_{n\ge 2,\Omega(n)=k} \left(1-\frac{1}{n^s}\right)
,
\quad
\zeta_1(s)=\zeta(s)
,
\quad \Re s>1.
\label{eq.zetaksDef}
\end{equation}
\end{defn}
With reference to (\ref{eq.P1}) there are logarithms, the logarithmic derivative,
\begin{equation}
\log \zeta_k(s) = \sum_{j=1}^\infty \frac{1}{j} P_k(js)
,\quad
\frac{ \zeta_k'(s)}{\zeta_k(s)} = \sum_{j=1}^\infty P_k'(js)
,
\end{equation}
and their M\"obius inversion
\begin{equation}
P_k(s) = \sum_{j=1}^\infty \frac{\mu(j)}{j} \log \zeta_k(js)
.
\label{eq.Pksofz2}
\end{equation}
This definition factorizes each of the constants of (\ref{eq.JVP}),
\begin{equation}
\prod_{n=2}^\infty \left(1-\frac{1}{n^s}\right)
=\prod_{k=1}^\infty \frac{1}{\zeta_k(s)}
.
\end{equation}
The unrestricted host product over all integers is evaluated with
(\ref{eq.ntoGamma}),
\begin{equation}
\prod_{n=2}^\infty \left(1-\frac{1}{n^s}\right)
=
\prod_{n=1}^\infty \frac{(n+1)^s-1}{(n+1)^s}
=
\frac{1}{\prod_{l=1}^{s-1}\Gamma(2-e^{2\pi il/s})}
\end{equation}
\begin{equation}
=\left\{
\begin{array}{ll}
1/2,& s=2\\
0.809396597366290\ldots,& s=3,\\
0.919019477593744\ldots,& s=4,\\
0.963256561757559\ldots,& s=5,\\
0.982684277742192\ldots,& s=6,\\
0.991654953472834\ldots,& s=7,\\
0.995923315077783\ldots,& s=8,\\
0.997991715347709\ldots,& s=9,\\
0.999005442480989\ldots,& s=10.\\
\end{array}
\right.
=\left\{
\begin{array}{ll}
\frac{\cosh(\pi\surd 3/2)}{3\pi},& s=3,\\
\frac{\sinh(\pi)}{4\pi},& s=4,\\
\frac{\cosh^2(\pi\surd 3/2)}{6\pi^2},& s=6,\\
\frac{\sinh \pi\left[(\sin \tau\cosh\tau)^2+(\cos\tau\sinh \tau)^2\right]}{8\pi^3}
,& s=8,
\end{array}
\right.
\label{eq.Zofs}
\end{equation}
where $\tau\equiv \pi/\surd 2$.
\begin{table}
\caption{$\zeta_k$ defined in (\ref{eq.zetaksDef}).
The first block is the familiar Riemann zeta function \cite[Table 23.2]{AS}.
As $s\to \infty$, $\zeta_k \sim 1+2^{-ks}$.
}
\begin{tabular}{|l|l|l|}
\hline
$k$ & $s$ & $\zeta_k(s)$ \\
\hline
1 &  2& 1.644934066848226436472415166646025189218949901206798437735
$\ldots$
\\
1 &  3& 1.202056903159594285399738161511449990764986292340498881792
$\ldots$
\\
1 &  4& 1.082323233711138191516003696541167902774750951918726907682
$\ldots$
\\
1 &  5& 1.036927755143369926331365486457034168057080919501912811974
$\ldots$
\\
1 &  6& 1.017343061984449139714517929790920527901817490032853561842
$\ldots$
\\
1 &  7& 1.008349277381922826839797549849796759599863560565238706417
$\ldots$
\\
1 &  8& 1.004077356197944339378685238508652465258960790649850020329
$\ldots$
\\
\hline
2 &  2& 1.154135429131192212753136476082653062021377019769166311601
$\ldots$
\\
2 &  3& 1.024230611826986151158175158755009852679023950490214554774
$\ldots$
\\
2 &  4& 1.005015172899917179827401698652291294164076627857960524772
$\ldots$
\\
2 &  5& 1.001137144097153269444733210594786272022830561556805839703
$\ldots$
\\
2 &  6& 1.000268773319796626045829352019049130925353448777275622774
$\ldots$
\\
2 &  7& 1.000064937119208710399247524597916121342642545376082696815
$\ldots$
\\
2 &  8& 1.000015888762704374180969521932691056219718849573814329451
$\ldots$
\\
2 &  9& 1.000003917599435801201288001083572332107127973809390534197
$\ldots$
\\
2 &  10& 1.000000970605519540457698447691588436106985664105867746091
$\ldots$
\\
2 &  11& 1.000000241217232784709472971654363631286102379240142530861
$\ldots$
\\
2 &  12& 1.000000060068608308578085912628099422424825142999323470824
$\ldots$
\\
\hline
3 &  2& 1.039432429030409444149806920521673410679212706746148280292
$\ldots$
\\
3 &  3& 1.003056125733692715934271390834372275666674738470640544373
$\ldots$
\\
3 &  4& 1.000314507983058618456827094397066630425050030824169155161
$\ldots$
\\
3 &  5& 1.000035578360205492419327661940592089489174843096220340325
$\ldots$
\\
3 &  6& 1.000004201291691977098551061424155872921133219822741524992
$\ldots$
\\
3 &  7& 1.000000507389042251328731405257671868967853054036047135925
$\ldots$
\\
3 &  8& 1.000000062068139946907813780148372010808995726576667040951
$\ldots$
\\
3 &  9& 1.000000007651664785056457877351300502168705881798866985821
$\ldots$
\\
3 &  10& 1.000000000947860936856417390158481965876008262486022523638
$\ldots$
\\
3 &  11& 1.000000000117782004432039022814899184474717980169845755406
$\ldots$
\\
3 &  12& 1.000000000014665193378639871887901001332108989045782221764
$\ldots$
\\
\hline
4 &  2& 1.010069659181975191078741060439035427876588103787130676067
$\ldots$
\\
4 &  3& 1.000384011077316192420561416488116578273502799964152453503
$\ldots$
\\
4 &  4& 1.000019679948504124660855296987073387779196556022576263542
$\ldots$
\\
4 &  5& 1.000001112106579607264579274143875169225449518083026268097
$\ldots$
\\
4 &  6& 1.000000065648673607782540354304127473976972447884634538215
$\ldots$
\\
4 &  7& 1.000000003964020108633384822995698964455986578036971179818
$\ldots$
\\
4 &  8& 1.000000000242454206776351387845438137782168399984862743746
$\ldots$
\\
4 &  9& 1.000000000014944664405752269172969167215376338782896827053
$\ldots$
\\
4 &  10& 1.000000000000925645528057890635386834145659663652892281940
$\ldots$
\\
4 &  11& 1.000000000000057510745364142948158605083885270190381157305
$\ldots$
\\
4 &  12& 1.000000000000003580369489719342307308508350305504446132847
$\ldots$
\\
\hline
\end{tabular}
\label{tab.zetak}
\end{table}

Cyclotomic factorizations (avoiding divergent expressions,  ie,
using only polynomial orders of 2 or higher) generalize (\ref{eq.zetadoubl})
if $x=1/n$ is set to visualize the structure:
\begin{itemize}
\item
From $1-x^4=(1-x^2)(1+x^2)$ we generate
\begin{equation}
\prod_{n\ge 2,\Omega(n)=k} \left(1-\frac{1}{n^4}\right) =
\prod_{n\ge 2,\Omega(n)=k} \left(1-\frac{1}{n^2}\right)
\left(1+\frac{1}{n^2}\right)
.
\end{equation}
Division leads to
\begin{equation}
\prod_{n\ge 2,\Omega(n)=k}
\left(1+\frac{1}{n^2}\right)
=\frac{\zeta_k(2)}{\zeta_k(4)}
,
\end{equation}
and in general to
\begin{equation}
\prod_{n\ge 2,\Omega(n)=k}
\left(1+\frac{1}{n^s}\right)
=\frac{\zeta_k(s)}{\zeta_k(2s)}
,
\quad
\therefore
\prod_{n\ge 2,\Omega(n)=k}
\frac{n^s-1}{n^s+1}
=\frac{\zeta_k(2s)}{\zeta_k^2(s)}
.
\label{eq.Zofzeta2}
\end{equation}
In that sense, the $\zeta_k$ keep up to the promise to generate the rational
polynomials of Section \ref{sec.intro}.
\item
From $1-x^6=(1-x^3)(1+x^3)=(1-x^2)(1+x^2+x^4)$ we generate
\begin{equation}
\prod_{n\ge 2,\Omega(n)=k}
\left(1+\frac{1}{n^2}+\frac{1}{n^4}\right)
=\frac{\zeta_k(2)}{\zeta_k(6)}
,
\end{equation}
and generally from $1-x^{2s}=(1-x^2)(1+x^2+x^4+\cdots +x^{2s-2})$ 
\begin{equation}
\prod_{n\ge 2,\Omega(n)=k}
\left(1+\sum_{j=1}^s \frac{1}{n^{2j}}\right)
=\frac{\zeta_k(2)}{\zeta_k(2s+2)}
.
\end{equation}

\item
The three different ways of grouping factors
on the right-hand side of $(1-x^8)=(1-x^2)(1+x^2)(1+x^4)$ establish in addition
\begin{equation}
\prod_{n\ge 2,\Omega(n)=k}
\left(1-\frac{1}{n^2}+\frac{1}{n^4}-\frac{1}{n^6}\right)
=\frac{\zeta_k(4)}{\zeta_k(2)\zeta_k(8)}
,
\end{equation}
a special case of
\begin{equation}
\prod_{n\ge 2,\Omega(n)=k}
\left(1+\sum_{j=1}^s \left(-\frac{1}{n^2}\right)^j\right)
=\left\{
\begin{array}{ll}
\displaystyle \frac{\zeta_k(4)}{\zeta_k(2) \zeta_k(2s+2)}
  ,& s\, \mathrm{odd}; \\
\displaystyle \frac{\zeta_k(4)\zeta_k(2s+2)}{\zeta_k(2) \zeta_k(4s+4)}
  ,& s\, \mathrm{even}.
\end{array}
\right.
\end{equation}
\item
From $1-x^9=(1-x^3)(1+x^3+x^6)$  by division through $1-x^3$,
\begin{equation}
\prod_{n\ge 2,\Omega(n)=k}
\left(1+\frac{1}{n^{3j}}+\frac{1}{n^{6j}}\right)
=\frac{\zeta_k(3j)}{\zeta_k(9j)}
.
\end{equation}
\end{itemize}
Other formats are products of these in disguise. The elementary
examples are
\begin{equation}
\prod_{n\ge 2,\Omega(n)=k}
\left(1-\frac{1\mp 2n^s}{n^{2s}}\right)
=
\prod_{n\ge 2,\Omega(n)=k}
\left(1\pm \frac{1}{n^s}\right)^2,
\end{equation}
\begin{equation}
\prod_{n\ge 2,\Omega(n)=k}
\left(1-\frac{n^{s+1}\pm n^s\mp 1}{n^{2s+1}}\right)
=
\prod_{n\ge 2,\Omega(n)=k}
\left(1- \frac{1}{n^s}\right)
\left(1\pm \frac{1}{n^{s+1}}\right),
\end{equation}
to be continued in Appendix \ref{sec.hybr}.

\section{Artin's Constant}

\begin{defn}
(Artin's Constants of order $r$) \cite{WrenchMCom15}
\begin{equation}
A^{(r)}\equiv \prod_{n=2}^\infty \left(1-\frac{1}{n^r(n-1)} \right)
=\prod_{k=1}^\infty A_k^{(r)};\quad 
A_k^{(r)}\equiv \prod_{\substack{n=2\\ \Omega(n)=k}}^\infty 
\left(1-\frac{1}{n^r(n-1)} \right)
\label{eq.ArtDef}
.
\end{equation}
\end{defn}
\begin{gather}
\log A^{(r)}
=
\sum_{n=2}^\infty \log \left(1-\frac{1}{n^r(n-1)}\right)
=
- \sum_{n=2}^\infty \sum_{s=1}^\infty \frac{1}{s n^{rs}(n-1)^s}
\label{eq.Ar}
\\
=
- \sum_{n=2}^\infty \sum_{s=1}^\infty \sum_{l=0}^\infty \binom{-s}{l}\frac{(-1)^l}{s n^{(r+1)s+l}}
=
- \sum_{s=1}^\infty \sum_{l=0}^\infty \frac{(s)_l}{sl!}\left[\zeta([r+1]s+l)-1\right]
\nonumber
\\
=
\sum_{s=2}^\infty \sum_{j=1}^{\lfloor s/(1+r)\rfloor} \frac{1}{j}\binom{s-rj-1}{j-1}[1-\zeta(s)]
.
\nonumber
\end{gather}
\begin{table}
\caption{Artin's constants defined in (\ref{eq.ArtDef}).
Where the $k$-column is empty, the value is $A^{(r)}$, else $A_k^{(r)}$.
}
\begin{tabular}{|l|l|l|}
\hline
$r$ & $k$ & $A^{(r)}$, $A_k^{(r)}$ \\
\hline
1 &   & 0.296675134743591034570155020219142864864831519178947890816
$\ldots$
\\
1 & 1 & 0.373955813619202288054728054346416415111629248606150042094
$\ldots$
\\
1 & 2 & 0.839042154274468600768462111194541254928307166760882733000
$\ldots$
\\
1 & 3 & 0.958752116435730927714740256578928612659490448502359901592
$\ldots$
\\
1 & 4 & 0.989628867166427665504322837457924308057557589350296534844
$\ldots$
\\
\hline
2 &   & 0.673917363376357541664408979322634438564759812312671736792
$\ldots$
\\
2 & 1 & 0.697501358496365903284670350820922924073153946214515395354
$\ldots$
\\
2 & 2 & 0.969932325001525316214920207789129575961145794796696088006
$\ldots$
\\
2 & 3 & 0.996598927480241273419159046329894692291010391011783820658
$\ldots$
\\
2 & 4 & 0.999595278586535535637452493248336453083650632412674049887
$\ldots$
\\
\hline
3 &   & 0.850670630791104353750309521250006234999150598195442830656
$\ldots$
\\
3 & 1 & 0.856540444853542174426167984135953882166572800317652140325
$\ldots$
\\
3 & 2 & 0.993521589710505460675409269241416429401115078677815660188
$\ldots$
\\
3 & 3 & 0.999645238332613367730206639120726777503960574831358345008
$\ldots$
\\
\hline
4 &   & 0.929838473954346852238318469534553548944908305482253635236
$\ldots$
\\
4 & 1 & 0.931265184160004334389237205550676982558423734587801059016
$\ldots$
\\
4 & 2 & 0.998509500607573754587150213578131277848209932996862159036
$\ldots$
\\
4 & 3 & 0.999959643476392507167505941218783984697629873351553896969
$\ldots$
\\
\hline
5 &   & 0.966321276366930291670339804179360258225974383645878751173
$\ldots$
\\
5 & 1 & 0.966668868596777512740328372930016264211423822118193979007
$\ldots$
\\
5 & 2 & 0.999645271339337636005026101684121966092871540313156977621
$\ldots$
\\
5 & 3 & 0.999995220535319372944786385639468021407166618867177761801
$\ldots$
\\
\hline
\end{tabular}
\label{tab.Art}
\end{table}
This demonstrates how
\begin{equation}
\log A_k^{(r)}
=
- \sum_{s=2}^\infty \sum_{j=1}^{\lfloor s/(1+r)\rfloor} \frac{1}{j}\binom{s-rj-1}{j-1}P_k(s)
\label{eq.logAkr}
\end{equation}
is derived; the actual numerical
evaluation of $A^{(r)}$ is done easier via (\ref{eq.ntoGamma}).
Occasionally these $\Gamma$-functions simplify:
\begin{equation}
A^{(1)}= -\frac{\sin(\pi\phi)}{\pi},\quad  \phi\equiv \frac{\sqrt{5}+1}{2}.
\end{equation}

\begin{defn}
(Decremented Generalized Lucas Sequences)
\begin{equation}
a_{r,s}\equiv
s\sum_{j=1}^{\lfloor s/(r+1)\rfloor} \frac{1}{j}\binom{s-jr-1}{j-1}
.
\label{eq.arsdef}
\end{equation}
The attribute ``decremented'' stresses that $1+a_{r,s}$ are the reference
values that one would find, for example, in the Online Encyclopedia
of Integer Sequences \cite{EIS}.
\end{defn}
These integer sequences provide the notational shortcut
\begin{equation}
\log A_k^{(r)}= -\sum_{s=2}^\infty \frac{a_{r,s}}{s}P_k(s),
\label{eq.AkrofPks}
\end{equation}
and have recurrences
\begin{equation}
a_{r,s}=0,\quad s\le 1;\quad
a_{r,s} = 
2 a_{r,s-1}-a_{r,s-2}+a_{r,s-r-1}-a_{r,s-r-2},\quad  r\ge 1
\label{eq.arsrec}
\end{equation}
and generating functions
\begin{equation}
\sum_{s=0}^\infty a_{r,s} x^s = \frac{x^{1+r}(1+r-rx)}{(1-x)(1-x-x^{1+r})}
=
-r-\frac{1}{1-x}+\frac{1+r-rx}{1-x-x^{1+r}}
.
\label{eq.arsgf}
\end{equation}

Values for $s\ge 2$ are
\begin{eqnarray*}
a_{1,s} &=& 2,3,6,10,17,28,46,75,122,198,321,520,842,1363,2206,3570,\ldots
\\
a_{2,s} &=& 0,3,4,5,9,14,20,30,45,66,97,143,210,308,452,663,972,1425,2089,
\ldots
\\
a_{3,s} &=& 0,0,4,5,6,7,12,18,25,33,46,65,91,125,172,238,330,456,629,868,
\ldots
\end{eqnarray*}

The argument of the product (\ref{eq.ArtDef}) admits an exponential product
expansion \cite{MoreeMM101,Niklasch}
\begin{equation}
1-\frac{1}{n^r(n-1)} = \prod_{j=1}^\infty \left(1-\frac{1}{n^j}\right)^{\gamma_{r,j}^{(A)}}
.
\label{eq.gammADef}
\end{equation}
From the Laurent expansion
\begin{equation}
\frac{1}{n^r(n-1)} 
=
\sum_{j=0}^\infty \frac{1}{n^{r+1+j}}
.
\end{equation}
we see that this is (up to a sign flip) the inverse Euler transform
of the all-1 sequence padded with a short string of initial zeros \cite{CameronJIS3,BernsteinLAA226}.

\begin{eqnarray}
\gamma_{1,j}^{(A)} &=& 0, 1, 1, 1, 2, 2, 4, 5, 8, 11, 18, 25, 40, 58, 90, 135, 210, 316, 492,\ldots
\\
\gamma_{2,j}^{(A)} &=& 0, 0, 1, 1, 1, 1, 2, 2, 3, 4, 6, 7, 11, 14, 20, 27, 39, 52, 75, 102, 145,\ldots
\\
\gamma_{3,j}^{(A)} &=& 0, 0, 0, 1, 1, 1, 1, 1, 2, 2, 3, 3, 5, 6, 8, 10, 14, 17, 24, 30, 41, 53,\ldots
\end{eqnarray}

Whereas (\ref{eq.logAkr}) expands $A_k^{(r)}$ in the $P_k$ basis,
(\ref{eq.gammADef}) generates them in $\zeta_k$ basis:
\begin{equation}
A_k^{(r)}
=
\prod_{\Omega(n)=k} \prod_j \left(1-\frac{1}{n^j}\right)^{\gamma_{r,j}^{(A)}}
=
\prod_j \zeta_k(j)^{-\gamma_{r,j}^{(A)}}
.
\end{equation}
The same result could be obtained by plugging (\ref{eq.Pksofz2}) into
the right hand side of (\ref{eq.AkrofPks}), then exponentiation,
which reveals the M\"obius pair
\begin{equation}
\frac{1}{j}\sum_{l\mid j}\mu(l) a_{r,j/l} = \gamma_{r,j}^{(A)};\quad
a_{r,s} = \sum_{l\mid s}l \gamma_{r,l}^{(A)}.
\label{eq.mobA}
\end{equation}
\begin{rem}
This connection between the two types of coefficients
does not depend on the mediation by $P$ nor on summation over $n$,
and has been observed before \cite{HuaJCTA79,LiJCTA115,deReyniaAcAr119}.
Supposed any $\gamma_j$ are defined in
\begin{equation}
X\equiv \prod_j\left(1-\frac{1}{n^j}\right)^{\gamma_j},
\end{equation}
which implies
\begin{equation}
\log X=\sum_j \log\left(1-\frac{1}{n^j}\right)^{\gamma_j}
= \sum_j \gamma_j\log\left(1-\frac{1}{n^j}\right)
= -\sum_j \gamma_j \sum_{s\ge 1} \frac{1}{jn^{sj}}.
\end{equation}
Let furthermore $g_k$ be defined via
\begin{equation}
\log X\equiv -\sum_k \frac{g_k}{k}\frac{1}{n^k},
\end{equation}
then---by comparison of coefficients of equal powers $k=sj$ of $n$---
\begin{equation}
\frac{g_k}{k}=\sum_{j\mid k} \frac{\gamma_j}{k/j}
\quad
\therefore
\quad
g_k=\sum_{j\mid k} j\gamma_j
\quad
\therefore
\quad
\gamma_j = \frac{1}{j}\sum_{l\mid j}\mu(l)g_{j/l}.
\label{eq.mob}
\end{equation}
\end{rem}

\section{Twin Prime Constants}
\begin{defn}
(Twin Prime Constants of order $r$)
\begin{equation}
T^{(r)}\equiv \prod_{n=3}^\infty \left(1-\frac{1}{(n-1)^r} \right)
=\prod_{k=1}^\infty T_k^{(r)};\quad 
T_k^{(r)}\equiv \prod_{\substack{n=3\\ \Omega(n)=k}}^\infty 
\left(1-\frac{1}{(n-1)^r} \right)
\label{eq.TprDef}
.
\end{equation}
\end{defn}
\begin{gather}
\log T^{(r)}
=
\sum_{n=3}^\infty \log \left(1-\frac{1}{(n-1)^r}\right)
=
- \sum_{n=3}^\infty \sum_{s=1}^\infty \frac{1}{s (n-1)^{rs}}
\label{eq.logTr}
\\
=
- \sum_{n=3}^\infty \sum_{s=1}^\infty \sum_{l=0}^\infty \binom{-rs}{l}
  \frac{(-1)^l}{s n^{rs+l}}
=
- \sum_{n=3}^\infty \sum_{s=1}^\infty \sum_{l=0}^\infty \frac{(rs)_l}{sl! n^{rs+l}}
\nonumber
\\
=
- \sum_{s=1}^\infty \sum_{l=0}^\infty
  \frac{\Gamma(rs+l)}{\Gamma(rs)sl!}\left[\zeta(rs+l)-1-\frac{1}{2^{rs+l}}\right]
\nonumber
\\
=
- \sum_{s=r}^\infty \sum_{j=1}^{\lfloor s/r\rfloor}
  \frac{1}{j}\binom{s-1}{rj-1}\left[\zeta(s)-1-\frac{1}{2^s}\right]
=
- \sum_{n=2}^\infty \sum_{s=1}^\infty \frac{1}{s n^{rs}}
=
\sum_{s=1}^\infty \frac{1}{s}\left[1-\zeta(rs)\right]
.
\nonumber
\end{gather}
This is a sum rule for the zeta function. The actual values of $T^{(r)}$
duplicate those of (\ref{eq.Zofs}) because
\begin{equation}
T^{(r)} = 
\prod_{n=2}^\infty (1-\frac{1}{n^r})
.
\end{equation}
Repeating (\ref{eq.logTr}) yields
\begin{equation}
\log T_k^{(r)}
=
- \sum_{s=1}^\infty \sum_{j=1} ^{\lfloor s/r\rfloor}\frac{1}{j}
 \binom{s-1}{rj-1}\times
\left\{ \begin{array}{ll} P_k(s), & k>1 \\
P(s)-\frac{1}{2^s}, & k=1 \\
\end{array}
\right.
\label{eq.Akr}
.
\end{equation}
One can rewrite this as
\begin{equation}
\log T_k^{(r)}
=
- \sum_{s=1}^\infty \frac{r}{s}t_{r,s}
\times
\left\{ \begin{array}{ll} P_k(s), & k>1 \\
P(s)-\frac{1}{2^s}, & k=1 \\
\end{array}
\right.
.
\end{equation}
by introducing integer sequences $t_{r,s}$ via
\begin{defn}
(Binomial transform of aerated all-1 sequences)
\begin{equation}
t_{r,s}\equiv
\sum_{j=1} ^{\lfloor s/r\rfloor}\frac{s}{rj}
 \binom{s-1}{rj-1}
.
\end{equation}
\end{defn}
Generating functions are
\begin{equation}
\sum_{s=1}^\infty t_{r,s} x^s = 
\left\{
\begin{array}{ll}
\frac{x^2}{(1-x)(1-2x)},& r=2,\\
\frac{x^3}{(1-x)(1-2x)(1-x+x^2)},& r=3, \\
\frac{x^4}{(1-x)(1-2x)(1-2x+2x^2)},& r=4, \\
\frac{x^5}{(1-x)(1-2x)(1-3x+4x^2-2x^3+x^4)},& r=5,\\
\frac{x^6}{(1-x)(1-2x)(1-4x+7x^2-6x^3+3x^4)},& r=6,\\
\frac{x^7}{(1-x)(1-2x)(1-5x+11x^2-13x^3+9x^4-3x^5+x^6)},& r=7.
\end{array}
\right.
\end{equation}
\begin{eqnarray}
t_{2,s}&=& 0,1,3,7,15,31,63,127,255,511,1023,2047,4095,8191,16383
\ldots
\\
t_{3,s}&=& 0,0,1,4,10,21,42,84,169,340,682,1365,2730,5460,10921
\ldots
\\
t_{4,s}&=& 0,0,0,1,5,15,35,71,135,255,495,991,2015,4095,8255,16511
\ldots
\end{eqnarray}
The exponents $\gamma_{r,j}^{(T)}$ of the
inverse Euler transformation are defined as
\begin{equation}
1-\frac{1}{(n-1)^r} = \prod_{j=1}^\infty \left(1-\frac{1}{n^j}\right)^{\gamma_{r,j}^{(T)}}
,
\quad
T_k^{(r)} = \prod_{j\ge 2} \zeta_k^{-\gamma_{r,j}^{(T)}}.
\end{equation}
Examples at indices $j\ge 2$ are
\begin{eqnarray}
\gamma_{2,j}^{(T)} &=& 1, 2, 3, 6, 9, 18, 30, 56, 99, 186, 335, 630, 1161, 2182, 4080,
\ldots ,
\\
\gamma_{3,j}^{(T)} &=& 0, 1, 3, 6, 10, 18, 30, 56, 99, 186, 335, 630, 1161, 2182, 4080,
\ldots ,
\\
\gamma_{4,j}^{(T)} &=& 0, 0, 1, 4, 10, 20, 35, 60, 100, 180, 325, 620, 1160, 2200, 4110,
\ldots ,
\\
\gamma_{5,j}^{(T)} &=& 0, 0, 0, 1, 5, 15, 35, 70, 126, 215, 355, 605, 1065, 2002, 3855,
\ldots .
\end{eqnarray}
$\gamma_{2,j}^{(T)}$ and $\gamma_{3,j}^{(T)}$ differ only
by one at $j=2$, $3$ and $6$, caused by
\begin{equation}
\frac{ T_k^{(2)} }{ T_k^{(3)}}
=
\frac{ \zeta_k(6)}{ \zeta_k(2)\zeta_k(3)}.
\end{equation}
Equivalent to (\ref{eq.mobA}) we have
\begin{equation}
\frac{r}{j}\sum_{l\mid j}\mu(l) t_{r,j/l} = \gamma_{r,j}^{(T)};\quad
t_{r,s} = \frac{1}{r}\sum_{l\mid s}l \gamma_{r,l}^{(T)}.
\end{equation}

\begin{table}
\caption{Constants defined in (\ref{eq.TprDef}).
}
\begin{tabular}{|l|l|l|}
\hline
$r$ & $k$ & $T_k^{(r)}$ \\
\hline
2 &  1& 0.6601618158468695739278121100145557784326233602847334133194
$\ldots$
\\
2 &  2& 0.8045082612474742003477609755804840258842916235384411373425
$\ldots$
\\
2 &  3& 0.9550572882298700493014462558041415734147334017474337181150
$\ldots$
\\
2 &  4& 0.9892051323193015395807816797711957195936925403036958416322
$\ldots$
\\
\hline
3 &  1& 0.8553921020033986085509619238173369142779617326343761891909
$\ldots$
\\
3 &  2& 0.9507543513576159108562851848128299370568198257735656138328
$\ldots$
\\
3 &  3& 0.9957472030934113787540687306667347661629225361751865362705
$\ldots$
\\
3 &  4& 0.9995497160906745509662751060657004870141623046341572587167
$\ldots$
\\
\hline
4 &  1& 0.9329472788050225154245474117724776691267802933315976922689
$\ldots$
\\
4 &  2& 0.9856013155080764416552558374557003044523758203959126811851
$\ldots$
\\
4 &  3& 0.9994884471180270660931453925494333667590005538061096053390
$\ldots$
\\
4 &  4& 0.9999751268767094622618062255772478610135549582924856688974
$\ldots$
\\
\hline
5 &  1& 0.9676641641449273928684588284554067020603735049395556363790
$\ldots$
\\
5 &  2& 0.9955134067603972588631219139726981275567920402258163509688
$\ldots$
\\
5 &  3& 0.999932947887271875019211361923619823854671560349775644808
$\ldots$
\\
5 &  4& 0.999998490529076495094430461010984203017854835139390019305
$\ldots$
\\
\hline
\end{tabular}
\label{tab.Twin}
\end{table}

\section{Quadratic Class Number}
A sign flip in a denominator of (\ref{eq.ArtDef}) provides
\begin{defn}
(Quadratic Class numbers of order $r$)
\begin{equation}
Q^{(r)}\equiv \prod_{n=2}^\infty \left(1-\frac{1}{n^r(n+1)} \right)
=\prod_{k=1}^\infty Q_k^{(r)};\quad 
Q_k^{(r)}\equiv \prod_{\substack{n=2\\ \Omega(n)=k}}^\infty 
\left(1-\frac{1}{n^r(n+1)} \right)
.
\label{eq.Qdef}
\end{equation}
\end{defn}
The special value
\begin{equation}
Q^{(1)}= -\frac{2\sin(\pi \phi)}{\pi}
= 2A^{(1)}
,\quad \phi\equiv \frac{\sqrt{5}+1}{2},
\end{equation}
is found with (\ref{eq.ntoGamma}).
Essentially duplicating the calculation in (\ref{eq.Ar}) we have
\begin{gather}
\log Q^{(r)}
=
\sum_{n=2}^\infty \log \left(1-\frac{1}{n^r(n+1)}\right)
=
- \sum_{n=2}^\infty \sum_{s=1}^\infty \frac{1}{s n^{rs}(n+1)^s}
\label{eq.Qr}
\\
=
\sum_{s=2}^\infty \sum_{j=1}^{\lfloor s/(1+r)\rfloor} \frac{(-1)^{s-(r+1)j}}{j}\binom{s-rj-1}{j-1}[1-\zeta(s)]
.
\nonumber
\end{gather}
\begin{table}
\caption{Constants defined in (\ref{eq.Qdef}).
Where the $k$-column is empty, the value is $Q^{(r)}$, else $Q_k^{(r)}$.
}
\begin{tabular}{|l|l|l|}
\hline
$r$ & $k$ & $Q^{(r)}$, $Q_k^{(r)}$ \\
\hline
1 &   & 0.593350269487182069140310040438285729729663038357895781633
\\
1 & 1 & 0.704442200999165592736603350326637210188586431417098049414
\ldots \\
1 & 2 & 0.884490615792645569156126530213936198197151790687002628832
\ldots \\
1 & 3 & 0.964758474366761979144138911837762037272106901540505442479
\ldots \\
1 & 4 & 0.990393442303116742704510324208970590429148773419205301803
\ldots \\
\hline
2 &   & 0.861465028009033072712078741634897482256486581138963373814
\ldots \\
2 & 1 & 0.881513839725170776928391822903227847129869257208076733670
\ldots \\
2 & 2 & 0.980376289243855939864938619579944255828568405323127482357
\ldots \\
2 & 3 & 0.997235390988435066639517992790322243014456604416487473356
\ldots \\
2 & 4 & 0.999634762421613451087385734159581784386619475630152927741
\ldots \\
\hline
3 &   & 0.943588586975055819450398681669287957384984288942155352763
\ldots \\
3 & 1 & 0.947733262143675375939521537654189613033631632317413852828
\ldots \\
3 & 2 & 0.995928002832231110330556534959384975519658530911655791754
\ldots \\
3 & 3 & 0.999717415940667569960015810241340983047451116331199660541
\ldots \\
3 & 4 & 0.999981370624803353143270623749614063667365657870395981306
\ldots \\
\hline
4 &   & 0.974894913359018345579696732396272984272954782066491224441
\ldots \\
4 & 1 & 0.975824153047668241679011436594799831971764971229212609442
\ldots \\
4 & 2 & 0.999080621322852565861019671107083386967236557968329389061
\ldots \\
4 & 3 & 0.999968172024482705270548857407825523042346422728579249674
\ldots \\
4 & 4 & 0.999998949807673405679237367569730617874479395755489593508
\ldots \\
\hline
5 &   & 0.988286601083665561883354705652451815692957697529740408374
\ldots \\
5 & 1 & 0.988504397741246908751106623851186664400958083275346188120
\ldots \\
5 & 2 & 0.999783481766640731283537327557402867111772830353712345808
\ldots \\
5 & 3 & 0.999996250854906273732628391333715607260086058728852979934
\ldots \\
5 & 4 & 0.999999938086969714053459138631507932333794604460092491194
\ldots \\
\hline
\end{tabular}
\label{tab.Quad}
\end{table}
Table \ref{tab.Quad} is calculated via
\begin{equation}
\log Q_k^{(r)}
=
- \sum_{s=2}^\infty \sum_{j=1}^{\lfloor s/(1+r)\rfloor} \frac{(-1)^{s-(r+1)j}}{j}\binom{s-rj-1}{j-1}P_k(s)
=
- \sum_{s=2}^\infty \frac{1}{s} q_{r,s} P_k(s)
.
\label{eq.logQkr}
\end{equation}
The previous line introduces
auxiliary integer sequences $q_{r,s}$
equivalent to (\ref{eq.arsdef}) with
\begin{defn}
(Binomial transforms of aerated alternating-1 sequences)
\begin{equation}
q_{r,s}\equiv
s\sum_{j=1}^{\lfloor s/(r+1)\rfloor} \frac{(-1)^{s-(r+1)j}}{j}\binom{s-jr-1}{j-1}
.
\end{equation}
\end{defn}
Lists at $s\ge 2$ are
\begin{eqnarray*}
q_{1,s} &=& 2, -3, 6, -10, 17, -28, 46, -75, 122, -198, 321, -520, 842, -1363, 2206, -3570
\ldots
\\
q_{2,s} &=& 0,3,-4,5,-3,0,4,-6,5,0,-7,13,-14,8,4,-17,24,-19,1,24,-44,46,-23,-20,
\ldots
\\
q_{3,s} &=& 0,0,4,-5,6,-7,12,-18,25,-33,46,-65,91,-125,172,-238,330,-456,629,
\ldots
\\
q_{4,s} &=& 0,0,0,5,-6,7,-8,9,-5,0,6,-13,21,-25,24,-17,3,19,-45,70,-88,92,-74,30,
\ldots
\end{eqnarray*}
Recurrences and generating function are sign-flipped variants of (\ref{eq.arsrec}) and (\ref{eq.arsgf}),
\begin{equation}
q_{r,s}=0,\quad s\le 1;\quad
q_{r,s} = 
-2 q_{r,s-1}-q_{r,s-2}+q_{r,s-r-1}+q_{r,s-r-2},\quad  r\ge 1
,
\end{equation}
\begin{equation}
\sum_{s=0}^\infty q_{r,s} x^s = \frac{x^{1+r}(1+r+rx)}{(1+x)(1+x-x^{1+r})}
=
-r-\frac{1}{1+x}+\frac{1+r+rx}{1+x-x^{1+r}}
.
\end{equation}
The zeta-expansion exponents $\gamma_{r,j}^{(Q)}$ are defined to satisfy
\begin{equation}
1-\frac{1}{n^r(1+n)} = \prod_{j=1}^\infty \left(1-\frac{1}{n^j}\right)^{\gamma_{r,j}^{(Q)}}
,\quad
Q_k^{(r)}=\prod_{j=2}^\infty \zeta_k^{-\gamma_{r,j}^{(Q)}}
.
\end{equation}
In the range $j\ge 2$ we find
\begin{eqnarray}
\gamma_{1,j}^{(Q)} =  1, -1, 1, -2, 3, -4, 5, -8, 13, -18, 25, -40, 62, -90, 135, -210,
\ldots
\\
\gamma_{2,j}^{(Q)} = 0, 1, -1, 1, -1, 0, 1, -1, 0, 0, 0, 1, -1, 0, 0, -1, 2, -1, 0, 1, -2, 2,
\ldots
\\
\gamma_{3,j}^{(Q)} = 0, 0, 1, -1, 1, -1, 1, -2, 3, -3, 3, -5, 7, -8, 10, -14, 19, -24, 30,
\ldots
\\
\gamma_{4,j}^{(Q)} = 0, 0, 0, 1, -1, 1, -1, 1, -1, 0, 1, -1, 1, -2, 2, -1, 0, 1, -2, 3, -4, 4.
\ldots
\end{eqnarray}
The analogue of (\ref{eq.mobA})
relates $q_{r,j}$ with $\gamma_{r,j}^{(Q)}$,
\begin{equation}
\frac{1}{j}\sum_{l\mid j}\mu(l) q_{r,j/l} = \gamma_{r,j}^{(Q)};\quad
q_{r,s} = \sum_{l\mid s}l \gamma_{r,l}^{(Q)}.
\end{equation}

\section{Feller-Tornier}
\begin{defn}
(Feller-Tornier Constants)
\begin{equation}
F^{(r)}\equiv \prod_{n=2}^\infty \left(1-\frac{2}{n^r} \right)
=\prod_{k=1}^\infty F_k^{(r)};\quad 
F_k^{(r)}\equiv \prod_{\substack{n=2\\ \Omega(n)=k}}^\infty 
\left(1-\frac{2}{n^r} \right)
.
\label{eq.FTdef}
\end{equation}
\end{defn}
Closed form expressions generated from the $n$ roots of $(n+1)^r-2$ via
(\ref{eq.ntoGamma}) are:
\begin{eqnarray}
F^{(2)} &=& - \frac{\sin(\pi\surd 2)}{\pi \surd 2},
\\
F^{(4)} &=& - \frac{\sinh(\pi \sqrt[4]{2})\sin(\pi\sqrt[4]{2})}{\pi^2 \surd 2}
,
\\
F^{(6)} &=& - \frac{\sin(\pi 2^{1/6})
\left[\cosh^2(\pi 2^{-5/6}\surd 3)-\cos^2(\pi 2^{-5/6})\right]
}{\pi^3\surd 2}
.
\end{eqnarray}
Table \ref{tab.FT} is produced accumulating
\begin{equation}
\log F_k^{(r)} =
\sum_{n\ge 2,\Omega(n)=k} \log\left(1-\frac{2}{n^r}\right)
=
-\sum_{j=1}^\infty \frac{2^j}{j}P_k(rj)
.
\label{eq.logFkr}
\end{equation}

\begin{table}
\caption{Constants defined in (\ref{eq.FTdef}).
Where the $k$-column is empty, the value is $F^{(r)}$, else $F_k^{(r)}$.
}
\begin{tabular}{|l|l|l|}
\hline
$r$ & $k$ & $F^{(r)}$, $F_k^{(r)}$ \\
\hline
2 &   & 0.216954294377476369356864039063437596599913299714241452950
\ldots \\
2 & 1 & 0.322634098939244670579531692548237066570950579665832709961
\ldots \\
2 & 2 & 0.746546589392028395045090198448531702307441861924359721333
\ldots \\
2 & 3 & 0.925267218004050491882112057709653605835184213438593596030
\ldots \\
2 & 4 & 0.980141423819685616602097797704982658698681266298832028699
\ldots \\
\hline
3 &   & 0.640575909221546133846815305854929804642156633120487523553
\ldots \\
3 & 1 & 0.676892737009881993610237326724389212797678397459788845273
\ldots \\
3 & 2 & 0.952981131498858970286460835642413457907291749293329766096
\ldots \\
3 & 3 & 0.993911463538908938115880607314708743157809479259098501054
\ldots \\
3 & 4 & 0.999232354384330965293278269072153569836124255417702962004
\ldots \\
\hline
4 &   & 0.840695833076274061650473710681177939252057653860271740567
\ldots \\
4 & 1 & 0.849732991384718766265053703629160439892820104242861046497
\ldots \\
4 & 2 & 0.990028758685057521928151579163269491391556722507312328948
\ldots \\
4 & 3 & 0.999371218596872044721227851298066869151465994816088393076
\ldots \\
4 & 4 & 0.999960641022410557951000962330383249199360087979036126331
\ldots
\\
\hline
5 &   & 0.926880857710656853256795504739069232524104014194213179034
\ldots \\
5 & 1 & 0.929059192959662815115245871984200623766376123420999266247
\ldots \\
5 & 2 & 0.997728614956085550822520064832476818851830182381404991433
\ldots \\
5 & 3 & 0.999928846129017559729606214395367084797186701036127260135
\ldots \\
5 & 4 & 0.999997775789625477894415002081234674890918717170135941506
\ldots \\
\hline
\end{tabular}
\label{tab.FT}
\end{table}
The zeta-expansion exponents $\gamma_{r,j}^{(F)}$ are
\begin{equation}
1-\frac{2}{n^r} = \prod_{j=1}^\infty \left(1-\frac{1}{n^j}\right)^{\gamma_{r,j}^{(F)}}
,
\quad
F_k^{(r)}
=
\prod_j \zeta_k(j)^{-\gamma_{r,j}^{(F)}}
\label{eq.FTgamma}
.
\end{equation}
(\ref{eq.logFkr}) in conjunction with (\ref{eq.Pksofz2}) reveals
\begin{equation}
\gamma_{r,j}^{(F)}=\left\{
\begin{array}{ll}
0 & , r \nmid j, \\
\displaystyle\frac{r}{j}\sum_{d\mid (j/r)} 2^d\mu(\frac{j}{rd})  & ,  r \mid j,
\end{array}
\right.
\end{equation}
so values in the range $j\ge 2$ start as
\begin{eqnarray}
\gamma_{2,j}^{(F)} = 2, 0, 1, 0, 2, 0, 3, 0, 6, 0, 9, 0, 18, 0, 30, 0, 56, 0, 99, 0, 186, 0, 335,
\ldots
\\
\gamma_{3,j}^{(F)} = 0, 2, 0, 0, 1, 0, 0, 2, 0, 0, 3, 0, 0, 6, 0, 0, 9, 0, 0, 18, 0, 0, 30, 0, 0, 56, 0,
\ldots
\end{eqnarray}
The others are obvious: the count of zero fillers---factors effectively
dropping out in (\ref{eq.FTgamma})---grows
simply as $r-1$.

\section{Hardy-Littlewood}

The Hardy-Littlewood constants, in the narrow sense, are the cases $k=1$ in
\begin{defn}\label{def.HL}
\begin{gather}
C^{(3)}\equiv \prod_{n=4}^\infty \left(1-\frac{3n-1}{(n-1)^3}\right)=\frac{2}{9} 
;\,
C_k^{(3)}\equiv \prod_{\substack{n\ge 4\\ \Omega(n)=k}} \left(1-\frac{3n-1}{(n-1)^3}\right)
;
\label{eq.HL3}
\\
C^{(4)}\equiv \prod_{n=5}^\infty \left(1-\frac{6n^2-4n+1}{(n-1)^4}\right)=\frac{3}{32}
;\,
C_k^{(4)}\equiv \prod_{\substack{n\ge 5\\ \Omega(n)=k}} \left(1-\frac{6n^2-4n+1}{(n-1)^4}\right)
;
\\
C^{(r)}\equiv \prod_{n>r} \frac{n^{r-1}(n-r)}{(n-1)^r}
=\frac{(r-1)!}{r^{r-1}}
;\,
C_k^{(r)}\equiv \prod_{\substack{n>r\\ \Omega(n)=k}}
\frac{n^{r-1}(n-r)}{(n-1)^r},\quad r\ge 3
\label{eq.HL5}
.
\end{gather}
\end{defn}
$T^{(2)}=C^{(2)}$ might be linked in here.
As before, we turn to the logarithm and define associate integer expansion
coefficients $c_s^{(r)}$ from their Laurent expansion,
\begin{gather}
\log C_k^{(r)}
=
\sum_{n> r,\Omega(n)=k}
\log\frac{n^{r-1}(n-r)}{(n-1)^r}
=
\sum_{n> r,\Omega(n)=k}
\log\frac{1-r/n}{(1-1/n)^r}
\\
=
\sum_{n> r,\Omega(n)=k}
\left[
\log\left(1-\frac{r}{n}\right)-r\log\left (1-\frac{1}{n}\right)
\right]
=
\sum_{n> r,\Omega(n)=k}
\left[
-\sum_{s\ge 1}\frac{r^s}{sn^s}+r\sum_{s\ge 1}\frac{1}{sn^s}
\right]
\nonumber
,
\end{gather}
which can be summarized after interchange of summations as
\begin{equation}
\log C_k^{(r)} =
-
\sum_{n> r,\Omega(n)=k}
\sum_{s=2}^\infty \frac{1}{s}\frac{c_s^{(r)}}{n^s},
\label{eq.logCkofcs}
\end{equation}
in a simple common format:
\begin{equation}
c^{(r)}_s\equiv \left\{
\begin{array}{ll}
0,& s<2; \\
r^s-r,& s\ge 2.
\end{array}
\right.
\end{equation}
Plugging (\ref{eq.P1}) into the right hand side of (\ref{eq.logCkofcs}),
we compensate for the fact that
the lower limits $n> r$ discard some $n$
for the prime ($k=1$) and semiprime ($k=2$) cases:
\begin{equation}
\log C_k^{(r)} =
- \sum_{s\ge 2} \frac{1}{s}c^{(r)}_{s}
\times\left\{ 
\begin{array}{ll}
\left[ P_k(s)-\frac{1}{2^s}-\frac{1}{3^s}-\frac{1}{5^s}\right], & r=5,6,\quad k=1  ; \\
\left[ P_k(s)-\frac{1}{2^s}-\frac{1}{3^s}\right], & r=3,4,\quad k=1;  \\
\left[ P_k(s)-\frac{1}{4^s}-\frac{1}{6^s}\right], & r=6,7,8,\quad k=2;   \\
\left[ P_k(s)-\frac{1}{4^s}\right], & r=4,5,\quad k=2;   \\
P_k(s), & r=3,\quad k=2;   \\
P_k(s), & 3\le r<8,\, k\ge 3. \\
\end{array}
\right.
\label{eq.logCkr}
\end{equation}
\begin{rem}
The convergence of these series is slow, given that $c_s^{(r)}$
grows roughly $\sim r^s$ and
the terms right from the brace
fall
roughly $\sim r^{-s}$.
The standard acceleration technique is to split the products (\ref{eq.HL5})
into $\prod_{n>r}=\prod_{r<n\le M}\cdot \prod_{n>M}$ with some free integer
$M$ of the order of some tens of $r$, to calculate the first of these
two products explicitly, and to subtract
all of the inverse powers of the $k$-almost primes below $M$
on the right hand side of (\ref{eq.logCkr}), so the
term to the right of the brace
falls off roughly $\propto M^{-s}$
\cite{SebahGourdon}.
\end{rem}

\begin{table}
\caption{
Hardy-Littlewood constants from 
equations (\ref{eq.HL3})--(\ref{eq.HL5}).
}
\begin{tabular}{|l|l|}
\hline
$k$ & $C_k^{(3)}$ \\
1 & 0.635166354604271207206696591272522417342065687332372450899
\ldots \\
2 & 0.424234558470737235218539671836177441479432573726566654172
\ldots \\
3 & 0.861978217115406397600389288178363010882226721530095958070
\ldots \\
4 & 0.967010333852598290029706098677367652117366812663687575395
\ldots \\
5 & 0.991986542483777613682589437104065646426247353044705005169
\ldots \\
\hline
    & $C_k^{(4)}$ \\
1 & 0.307494878758327093123354486071076853022178519950663928298
\ldots \\
2 & 0.461691758364773730232305524418356233105041873484187592372
\ldots \\
3 & 0.723165327592227885742644081506537901793428021559355399358
\ldots \\
4 & 0.933085922756286271428677500179124619333526568552797663388
\ldots \\
5 & 0.983814785274894677909254622541376985294671241029387304408
\ldots \\
\hline
    & $C_k^{(5)}$ \\
1 & 0.409874885088236474478781212337955277896358013254945469826
\ldots \\
2 & 0.199805231972458892888828284513805175888486651003235386176
\ldots \\
3 & 0.547976628430836989696044385054920027204371593584787789976
\ldots \\
4 & 0.887429384542023666239166084035101114754091610659986923636
\ldots \\
5 & 0.972787328073924092604485187403768369490797977225729870114
\ldots \\
\hline
    & $C_k^{(6)}$ \\
1 & 0.186614297358358396656924847944188337840073944945589304871
\ldots \\
2 & 0.298042020487754531592316128677284826210605850852999108669
\ldots \\
3 & 0.353138894039211423074594163633660113968205949561584457073
\ldots \\
4 & 0.830410751660277094955031322533872844216167956532800549112
\ldots \\
5 & 0.958867249262078290883709484646892100264447229912351544026
\ldots \\
\hline
\end{tabular}
\end{table}
Exponents $\gamma_{r,j}^{(C)}$ are defined
with the aim to decompose (\ref{eq.HL5}):
\begin{equation}
\frac{n^{r-1}(n-r)}{(n-1)^r}
\equiv \prod_{j\ge 1}
\left(1-\frac{1}{n^j}\right)^{\gamma_{r,j}^{(C)}}
,
\label{eq.gammCdef}
\end{equation}
related to $c_s^{(r)}$ via a M\"obius
transform as described in (\ref{eq.mob}) \cite{Niklasch}:
\begin{gather*}
\gamma_{3,j} = 3, 8, 18, 48, 116, 312, 810, 2184, 5880, 16104, 44220, 122640, 341484,\ldots,
\\
\gamma_{4,j} = 6, 20, 60, 204, 670, 2340, 8160, 29120, 104754, 381300, 1397740, 5162220,\ldots ,
\\
\gamma_{5,j} = 10, 40, 150, 624, 2580, 11160, 48750, 217000, 976248, 4438920, 20343700,\ldots
\end{gather*}
Eq.\ (\ref{eq.gammCdef}) rephrases the constants (\ref{eq.HL3})--(\ref{eq.HL5})
as
\begin{equation}
C_k^{(r)}=\prod_{j\ge 2} \lambda_k(r,j)^{-\gamma_{r,j}^{(C)}},
\end{equation}
where $\lambda_k$
takes into account that the products may have been defined
without the first one to three primes or semiprimes,
\begin{equation}
\lambda_k(r,j)=
\prod\limits_{\substack{n> r\\ \Omega(n)=k}}
\displaystyle \frac{1}{ (1-1/n^j)}
=
\zeta_k(j)\times \left\{
\begin{array}{ll}
(1-\frac{1}{2^j})(1-\frac{1}{3^j})(1-\frac{1}{5^j}), & r=5,\, k=1  ; \\
(1-\frac{1}{2^j})(1-\frac{1}{3^j}), & r=3,4,\, k=1;  \\
(1-\frac{1}{4^j}), & r=4,5,\, k=2;   \\
1, & r=3,\, k=2;   \\
1, & r=3,4,5,\, k\ge 3. \\
\end{array}
\right.
\end{equation}

\section{Summary}
The familiar Hardy-Littlewood, Artin's, Feller-Tornier and similar
constants are infinite products over the prime numbers. We have
generalized these to products over $k$-almost primes, and provide
tables
for low ranks and small
$k$.
The products over all values of $k$ are infinite ``host'' products,
easily evaluated as multi-gamma functions associated
with the roots of the defining rational polynomial.

\appendix
\section{Hybrids}\label{sec.hybr}
Hybrids are products and ratios of the constants discussed above, which
provide access to other forms of products. One example of such a
reduction is, see  (\ref{eq.Zofzeta2}),
\begin{equation}
\prod_{n\ge 2,\Omega(n)=k} \left( 1-\frac{3n^s+2}{n^{3s}}\right)
=
\prod_{n\ge 2,\Omega(n)=k} \left( 1-\frac{2}{n^s}\right)
\left( 1+\frac{1}{n^s}\right)^2
=
F_k^{(s)}\left( \frac{\zeta_k(s)}{\zeta_k(2s)}\right)^2.
\end{equation}
An---obviously incomplete---sample of these is:
\begin{equation}
\prod_{n\ge 2,\Omega(n)=k}
\left( 1+\frac{2n^s-1}{(n^s-1)^2}\right)
=
\zeta_k^2(s),
\end{equation}
\begin{equation}
\prod_{n\ge 2,\Omega(n)=k}
\left( 1+\frac{n^{s+l}+n^s-1}{n^{2s+l}-n^{s+l}-n^s+1}\right)
=
\zeta_k(s+l)\zeta_k(s),
\end{equation}
\begin{equation}
\prod_{n\ge 2,\Omega(n)=k}
\left( 1-\frac{1}{n^s+n^{s-1}+\cdots +n^2+n+1}\right)
=
\frac{\zeta_k(s+1)}{\zeta_k(s)},
\end{equation}
\begin{equation}
\prod_{n\ge 2,\Omega(n)=k}
\left( 1-\frac{n^l-1}{n^{s+l}-1}\right)
=
\frac{\zeta_k(s+l)}{\zeta_k(s)},
\end{equation}
\begin{equation}
\prod_{n\ge 2,\Omega(n)=k}
\left( 1-\frac{2n^s-1}{n^{2s}}\right)
=
\frac{1}{\zeta_k^2(s)},
\;\quad
\prod_{n\ge 2,\Omega(n)=k}
\left( 1-\frac{n^{s+l}+n^s-1}{n^{2s+l}}\right)
=
\frac{1}{\zeta_k(s+l)\zeta_k(s)},
\end{equation}
\begin{equation}
\prod_{n\ge 2,\Omega(n)=k}
\left( 1+\frac{1}{n(n^{s-1}+n^{s-2}+\cdots +n+1)}\right)
=
\frac{\zeta_k(s)}{\zeta_k(s+1)},
\end{equation}
\begin{equation}
\prod_{n\ge 2,\Omega(n)=k}
\left( 1-\frac{1}{n^s+1}\right)
=
\frac{\zeta_k(2s)}{\zeta_k(s)},
\end{equation}
\begin{equation}
\prod_{n\ge 2,\Omega(n)=k}
\left( 1-\frac{n^{(l-1)s}-1}{n^{ls}-1}\right)
=
\frac{\zeta_k(ls)}{\zeta_k(s)},
\end{equation}
\begin{equation}
\prod_{n\ge 2,\Omega(n)=k}
\left( 1+\frac{1}{n^s(n-1)-1}\right)
=
\frac{1}{ A_k^{(s)}},
\end{equation}
\begin{equation}
\prod_{n\ge 2,\Omega(n)=k}
\left( 1-\frac{1}{n^{s+2}-n^{s+1}-n+1}\right)
=
\zeta_k(s+1) A_k^{(s)},
\end{equation}
\begin{equation}
\prod_{n\ge 2,\Omega(n)=k}
\left( 1+\frac{2n+1}{n^3-2n-1}\right)
=
\frac{\zeta_k(2)}{A_k^{(1)}},
\end{equation}
\begin{equation}
\prod_{n\ge 2,\Omega(n)=k}
\left( 1+\frac{1}{n(n^{s+1}-n^s-1)}\right)
=
\frac{1}{\zeta_k(s+1)A_k^{(s)}},
\end{equation}
\begin{equation}
\prod_{n\ge 2,\Omega(n)=k}
\left( 1-\frac{2n+1}{n^3}\right)
=
\frac{ A_k^{(1)} } { \zeta_k(2) },
\quad
\prod_{n\ge 2,\Omega(n)=k}
\left( 1-\frac{2n^s+n^{s-1}+n^{s-2}+\cdots+n+1}{n^{2s+1}}\right)
=
\frac{ A_k^{(s)} } { \zeta_k(s+1) },
\end{equation}
\begin{equation}
\prod_{n\ge 2,\Omega(n)=k}
\left( 1+\frac{n-2}{n^{s+1}-n^s-n+1}\right)
=
\zeta_k(s)A_k^{(s)},
\end{equation}
\begin{equation}
\prod_{n\ge 2,\Omega(n)=k}
\left( 1-\frac{n-2}{n^{s+1}-n^s-1}\right)
=
\frac{1}{\zeta_k(s)A_k^{(s)} },
\end{equation}
\begin{equation}
\prod_{n\ge 2,\Omega(n)=k}
\left( 1-\frac{n^s+n^{s-1}+\cdots+n+1}{n^{2s}}\right)
=
\frac{ A_k^{(s)} } { \zeta_k(s) },
\end{equation}
\begin{equation}
\prod_{n\ge 2,\Omega(n)=k}
\left( 1-\frac{n^{s-1}+n^{s-2}+\cdots+n^2+2n+1}{n^{s+1}}\right)
=
\frac{ A_k^{(1)} } { \zeta_k(s) },
\end{equation}
\begin{equation}
\prod_{n\ge 2,\Omega(n)=k}
\left( 1+\frac{1}{n^s(n+1)-1}\right)
=
\frac{1}{ Q_k^{(s)}},
\end{equation}
\begin{equation}
\prod_{n\ge 2,\Omega(n)=k}
\left( 1+\frac{1}{n^{s+2}+n^{s+1}-n-1}\right)
=
\zeta_k(s+1) Q_k^{(s)},
\end{equation}
\begin{equation}
\prod_{n\ge 2,\Omega(n)=k}
\left( 1+\frac{2n+1}{n^3-2n+1}\right)
=
\frac{\zeta_k(2)}{Q_k^{(1)}},
\end{equation}
\begin{equation}
\prod_{n\ge 2,\Omega(n)=k}
\left( 1-\frac{1}{n(n^{s+1}+n^s-1}\right)
=
\frac{1}{\zeta_k(s+1)Q_k^{(s)}},
\end{equation}
\begin{equation}
\prod_{n\ge 2,\Omega(n)=k}
\left( 1-\frac{2n-1}{n^3}\right)
=
\frac{ Q_k^{(1)} } { \zeta_k(2) },
\quad
\prod_{n\ge 2,\Omega(n)=k}
\left( 1-\frac{2n^{s+1}+n^s-1}{n^{2s+1}(n+1)}\right)
=
\frac{ Q_k^{(s)} } { \zeta_k(s+1) },
\end{equation}
\begin{equation}
\prod_{n\ge 2,\Omega(n)=k}
\left( 1+\frac{n}{n^{s+1}+n^s-n-1}\right)
=
\zeta_k(s)Q_k^{(s)},
\end{equation}
\begin{equation}
\prod_{n\ge 2,\Omega(n)=k}
\left( 1-\frac{n}{n^{s+1}+n^s-1}\right)
=
\frac{1}{\zeta_k(s)Q_k^{(s)}},
\end{equation}
\begin{equation}
\prod_{n\ge 2,\Omega(n)=k}
\left( 1-\frac{n^2+n-1}{n^4}\right)
=
\frac{ Q_k^{(2)} } { \zeta_k(2) },
\quad
\prod_{n\ge 2,\Omega(n)=k}
\left( 1-\frac{n^{s+1}+2n^s-1}{n^{2s}(n+1)}\right)
=
\frac{ Q_k^{(s)} } { \zeta_k(s) },
\end{equation}
\begin{equation}
\prod_{n\ge 2,\Omega(n)=k}
\left( 1+\frac{2}{n^s-2}\right)
=
\frac{1}{F_k^{(s)}},
\end{equation}
\begin{equation}
\prod_{n\ge 2,\Omega(n)=k}
\left( 1-\frac{1}{n^s-1}\right)
=
\zeta_k(s)F_k^{(s)},
\quad
\prod_{n\ge 2,\Omega(n)=k}
\left( 1-\frac{2n^l-1}{n^{s+l}-1}\right)
=
\zeta_k(s+l)F_k^{(s)},
\end{equation}
\begin{equation}
\prod_{n\ge 2,\Omega(n)=k}
\left( 1+\frac{n^l-2}{n^l(n^s-1)}\right)
=
\zeta_k(s)F_k^{(s+l)},
\end{equation}
\begin{equation}
\prod_{n\ge 2,\Omega(n)=k}
\left( 1+\frac{3n^s-2}{n^{2s}-3n^s+2}\right)
=
\frac{\zeta_k(s)}{F_k^{(s)}},
\end{equation}
\begin{equation}
\prod_{n\ge 2,\Omega(n)=k}
\left( 1+\frac{1}{n^s-2}\right)
=
\frac{1}{\zeta_k(s)F_k^{(s)}},
\quad
\prod_{n\ge 2,\Omega(n)=k}
\left( 1+\frac{2n^l-1}{n^l(n^s-2)}\right)
=
\frac{1}{\zeta_k(s+l)F_k^{(s)}},
\end{equation}
\begin{equation}
\prod_{n\ge 2,\Omega(n)=k}
\left( 1-\frac{3n^s-2}{n^{2s}}\right)
=
\frac{ F_k^{(s)} }{ \zeta_k(s)}
,
\end{equation}
\begin{equation}
\prod_{n\ge 2,\Omega(n)=k}
\left( 1+\frac{n^l-1}{n^l(n^{s+1}-n^s-1)}\right)
=
\frac{ A_k^{(s+l)} }{ A_k^{(s)}}
,
\end{equation}
\begin{equation}
\prod_{n\ge 2,\Omega(n)=k}
\left( 1-\frac{n^l-1}{n^{s+l+1}-n^{s+l}-1}\right)
=
\frac{ A_k^{(s)} }{ A_k^{(s+l)}}
,
\end{equation}
\begin{equation}
\prod_{n\ge 2,\Omega(n)=k}
\left( 1-\frac{2n^{s+1}-1}{n^{2s}(n^2-1)}\right)
=
A_k^{(s)} Q_k^{(s)}
,
\end{equation}
\begin{equation}
\prod_{n\ge 2,\Omega(n)=k}
\left( 1-\frac{2}{n^{s+2}-n^s-n+1}\right)
=
\frac{A_k^{(s)}}{ Q_k^{(s)}}
,
\end{equation}
\begin{equation}
\prod_{n\ge 2,\Omega(n)=k}
\left( 1+\frac{2}{n^{s+2}-n^s-n-1}\right)
=
\frac{Q_k^{(s)}}{ A_k^{(s)}}
,
\end{equation}
\begin{equation}
\prod_{n\ge 2,\Omega(n)=k}
\left( 1+\frac{2n^{s+1}-n^s-2)}{n^{2s}(n-1)}\right)
=
A_k^{(s)} F_k^{(s)}
,
\end{equation}
\begin{equation}
\prod_{n\ge 2,\Omega(n)=k}
\left( 1-\frac{3n^{s+1}-2n^s-2)}{n^{2s+1}(n-1)}\right)
=
A_k^{(s)} F_k^{(s+1)}
,
\end{equation}
\begin{equation}
\prod_{n\ge 2,\Omega(n)=k}
\left( 1+\frac{2n-3}{n^{s+1}-n^s-2n+2}\right)
=
\frac{A_k^{(s)}}{ F_k^{(s)}}
,
\end{equation}
\begin{equation}
\prod_{n\ge 2,\Omega(n)=k}
\left( 1+\frac{n+2}{n^{s+2}+n^{s+1}-2n-2}\right)
=
\frac{Q_k^{(s)}}{ F_k^{(s+1)}}
,
\end{equation}
\begin{equation}
\prod_{n\ge 2,\Omega(n)=k}
\left( 1-\frac{n+2}{n(n^{s+1}+n^s-1)}\right)
=
\frac{F_k^{(s+1)}}{ Q_k^{(s)}}
,
\end{equation}
\begin{equation}
\prod_{n\ge 2,\Omega(n)=k}
\left( 1-\frac{2n^l+2}{n^{s+l}-2}\right)
=
\frac{F_k^{(s)}}{ F_k^{(s+l)}}
.
\end{equation}
Examples which involve $C_k^{(r)}$ have been left out for aesthetic
reasons, as
the dependence of the lower limit $n$ on $r$
in (\ref{eq.HL5}) leads to convoluted
subcase notation, like
\begin{equation}
\prod_{n\ge 5,\Omega(n)=k}
\left( 1+\frac{3}{n(n-4)}\right)
=
\left\{
\begin{array}{ll}
\frac{(27/16) C_k^{(3)}}{ C_k^{(4)}}, & k=2,\\
\frac{C_k^{(3)}}{ C_k^{(4)}}, & k\neq 2,
\end{array}
\right.
\end{equation}
which notices that the semiprime $n=4$ contributes to $C_k^{(3)}$
with a factor $16/27$
if $k=2$, but not otherwise.

\bibliographystyle{amsplain}
\bibliography{all}

\end{document}